# POISSON INVERSE PROBLEMS[1]

BY ANESTIS ANTONIADIS AND JÉREMIE BIGOT

*University Joseph Fourier and University Paul Sabatier*

In this paper we focus on nonparametric estimators in inverse problems for Poisson processes involving the use of wavelet decompositions. Adopting an adaptive wavelet Galerkin discretization, we find that our method combines the well-known theoretical advantages of wavelet–vaguelette decompositions for inverse problems in terms of optimally adapting to the unknown smoothness of the solution, together with the remarkably simple closed-form expressions of Galerkin inversion methods. Adapting the results of Barron and Sheu [*Ann. Statist.* **19** (1991) 1347–1369] to the context of log-intensity functions approximated by wavelet series with the use of the Kullback–Leibler distance between two point processes, we also present an asymptotic analysis of convergence rates that justifies our approach. In order to shed some light on the theoretical results obtained and to examine the accuracy of our estimates in finite samples, we illustrate our method by the analysis of some simulated examples.

**1. Introduction.** In this article the problem of estimating nonparametrically the intensity function of an indirectly observed nonhomogeneous Poisson process is considered. Such a problem arises when data (counts) are collected according to a Poisson process whose underlying intensity is indirectly related by a linear operator $K$ to the intensity (the object that we wish to estimate) of another Poisson process. This kind of indirect problem is referred to as a *Poisson inverse problem*. In rigorous probabilistic terms, let $F$ be a nonobservable Poisson process on a measure space $(E_0, \mathcal{B}(E_0), \mu_0)$ and let $tf(x)$ be its intensity function with respect to the measure $\mu_0$ on $\mathcal{B}(E_0)$; that is, for any set $A \in \mathcal{B}(E_0)$, the number of points of $F$ lying

Received September 2004; revised January 2006.

[1]Supported by funds from the IAP Research Network nr P5/24 of the Belgian Government (Federal Office for Scientific, Technical and Cultural Affairs) and the EC-HPRN-CT-2002-00286 Breaking Complexity Network.

AMS 2000 *subject classifications.* Primary 62G07; secondary 65J10.

in $A$ is a random variable $F(A)$ which is Poisson distributed with parameter $\int_A tf(x)\,d\mu_0(x)$, and for any finite family of disjoint measurable sets $A_1, \ldots, A_n$ of $E_0$, $F(A_1), \ldots, F(A_n)$ are independent random variables. For studying asymptotic properties, the function $f$, referred to as the scaled intensity function, is held fixed and the positive real $t$, referred to as the "observation time," increases. The observable data form another Poisson process $G$ on, possibly, another measure space $(E_1, \mathcal{B}(E_1), \mu_1)$ with an intensity function $th(y)$ with respect to a measure $\mu_1$. The scaled intensity functions $f$ and $h$, considered as elements of the separable Hilbert spaces $L^2(E_0, \mu_0)$ and $L^2(E_1, \mu_1)$, are related by an operator equation $h = Kf$ for some linear compact operator $K$ mapping $L^2(E_0, \mu_0)$ into $L^2(E_1, \mu_1)$. Observing the point process $G$ must be understood in a measure sense. Assuming therefore that for any $v \in L^2(E_1, \mu_1)$ we observe $\int v\,dG$, the natural goal is to estimate the scaled intensity $f$. In many applications, $K$ is an integral operator with a kernel representing the response of a measuring device; in the special case where this linear device is translation-invariant, $K$ reduces to a convolution operator. Examples range from all kinds of image deblurring models, mathematical models for positron emission tomography and nuclear magnetic resonance, or unfolding problems in stereology and high-energy physics, to cite only a few. Solving such problems, that is, recovering $f$, is often difficult since in cases which are of most interest scientifically, $K$ is not invertible; that is, $K^{-1}$ does not exist as a bounded linear operator so that a small perturbation in the data may lead to very different solutions to the recovery problem.

Related problems of inverse estimation for linear inverse problems with additive normal noise have been proposed in the literature, including smoothing kernel methods [12], smoothing spline methods [21, 22], Gauss–Chebyshev-type quadrature methods for solving integral equations [19] and singular value decomposition (SVD) methods [11, 23, 25], to cite only a few. Wavelet and multiscale analysis regularization methods for inverse problems have also recently received considerable attention in the statistics literature, exploiting the fact that wavelets provide unconditional bases for a large variety of smoothness spaces. Fan and Koo [9] have focused on nonparametric deconvolution density estimation based on wavelet techniques. Donoho [7] proposed the wavelet–vaguelette decomposition (WVD), which works by expanding the function $f$ in a wavelet series $\sum \langle f, \psi_{j,k} \rangle \psi_{j,k}$, constructing a corresponding vaguelette series for $Kf$ and then estimating the coefficients using a suitable thresholding approach. Donoho showed that a WVD is optimal in a minimax sense among all linear and nonlinear estimators for inverting certain types of linear operators, including the Radon transform. Kolaczyk [15] has numerically investigated the use of a WVD for tomographic reconstruction, whereas Abramovich and Silverman [1] have theoretically and numerically studied variants of the WVD. A drawback of these methods is that they are limited to special types of operators $K$ (essentially homogeneous operators or convolution-type operators under some additional technical assumptions) because one essentially needs to calculate precisely the $K^{-1}\psi_{j,k}$. Cohen, Hoffmann and Reiss [5] have explored the application of Galerkin-type methods to white-noise embedded inverse problems, using an appropriate but fixed wavelet basis. The underlying intuition is that the inversion process required by WVD methods needs only to be accurate to a certain error level if the object to be recovered is mostly smooth with some singularities, and therefore the inversion can be performed approximately using a Galerkin scheme. However, most of the techniques developed to date have been designed for Gaussian noise models and are not directly applicable in Poisson inverse problems.

For Poisson inverse problems an alternative approach has been proposed by Szkutnik [26] using a quasi-maximum likelihood (QML) histogram sieve estimator when restricting $h$ to step functions. Another recent attempt for solving related Poisson discrete inverse problems is a Bayesian multiscale framework for Poisson inverse problems proposed by Nowak and Kolaczyk [20], extending their earlier work for problems involving direct Poisson observations (see, e.g., [14, 16]) and based on a multiscale factorization of the Poisson likelihood function induced by recursive partitioning of the data space. Regularization of the solution is accomplished through usage of formal prior probability distributions in a Bayesian paradigm and the solution is a maximum a posteriori estimator, computed using the expectation–maximization (EM) algorithm. However, the inverse problems addressed by the above authors are *discrete* inverse problems (Poisson sampling from a discretized intensity related to a discretized version of the intensity of interest through multiplication by a matrix of transition probabilities), and the question up to which accuracy should the operators be discretized is not discussed. Similarly, the work of Cavalier and Koo [3] on hard threshold estimators in the tomographic data framework has shown that for a particular operator (the Radon transform) an extension of WVD methods for Poisson data is theoretically feasible. It is, however, worthwhile pointing out that the authors do not provide any computational algorithm for computing the estimate and do not address the problem of imposing positivity of the estimator since Poisson intensity functions are nonnegative by definition.

Encouraged by the developments cited above and inspired by the WVD methods for solving inverse problems, we explore in the sequel an alternative approach via wavelet-based decompositions combined with thresholding strategies that address adaptivity issues. Specifically, our framework extends the wavelet-Galerkin methods of Cohen, Hoffmann and Reiss [5] to the Poisson setting. Rates of convergence are derived for linear and nonlinear estimators, in analogy to classical wavelet estimators based on projections and thresholding, respectively. In order to ensure the positivity of the estimated intensity, the log-intensity is expanded in a wavelet basis. The derivation of our results takes place within an extension of the paradigm developed by Barron and Sheu [2] and involves the adaptation of recent techniques on concentration inequalities for suprema of integral functionals of Poisson processes which are analogous to Talagrand's inequalities for empirical processes. Although there are close similarities between wavelet-Galerkin techniques and earlier techniques based on WVD or VWD systems, the use of the wavelet-Galerkin machinery allows us to address inversion under a broad class of operators (i.e., not just homogeneous operators) and to take advantage of certain computational efficiencies.

The rest of this paper is organized as follows. While the Galerkin approach of Cohen, Hoffmann and Reiss [5] is relatively easy to describe when the inversion problem is a white-noise embedded problem, this is not the case for Poisson inverse problems. After fixing the notation and recalling some basic definitions, Section 2 contains an equivalent formulation of the wavelet-Galerkin approach for log-intensities that involves a notion of information projection similar to the one used by Barron and Sheu [2] for estimating a density, which is developed in Section 3. In Section 4 a linear estimator for the linear inverse problem at hand is proposed using the appropriate type of wavelets adapted to our case. In the spirit of wavelet denoising methods [1, 7], and in order to gain in adaptivity, we then improve, in Section 5, the estimator by applying a soft-threshold nonlinearity to the Galerkin-vaguelette coefficients. The last section is devoted to the numerical implementation of our procedures. We present the results of a small Monte Carlo experiment designed to study the finite-sample behavior of our estimates. Technical proofs are given in the Appendix.

## 2. Preliminaries and notation.

In this section we establish the notation and the general framework of the models which are adopted in this paper for the Poisson inverse problem formulated in the Introduction. For this purpose let $F$ and $G$ be two Poisson point processes on Borel measurable spaces $(E_0, \mathcal{B}(E_0), \mu_0)$ and $(E_1, \mathcal{B}(E_1), \mu_1)$, respectively. Associated with these Poisson processes are the *intensity measures* defined by

(a) $\lambda_F(B) = \mathbb{E}(F(B)) = \int_B tf(x)\, d\mu_0(x), B \in \mathcal{B}(E_0)$,

(b) $\lambda_G(B') = \mathbb{E}(G(B')) = \int_{B'} th(x)\, d\mu_1(x), B' \in \mathcal{B}(E_1)$,

where $t$ is an "observation time" which will tend to infinity in our asymptotic considerations. Observing the process $G$, we consider the problem of estimating the scaled intensity function $f$, when the scaled intensity $h$ results from the action of a compact self-adjoint positive definite operator $K: L^2(E_0, \mu_0) \to L^2(E_1, \mu_1)$ on the intensity function of the process $F$, that is, $h = Kf$. To simplify the notation, we will assume in the following without any loss of generality that the observation and unknown domains $E_0$ and $E_1$ are identical Borel subsets of $\mathbb{R}^d$ ($d \geq 1$), say $E$, and that $\mu_0 = \mu_1 = \mu$ where $\mu$ denotes Lebesgue measure. A discussion of how one can handle the case $E_0 \neq E_1$ or $K$ not self-adjoint positive definite is deferred to the end of this paper.

In order to estimate the unknown intensity function $f$, we will approximate the logarithm of the intensity by a standard wavelet basis function expansion. A notable advantage of using such an exponential family intensity estimation is that it forces positivity of the resulting estimator, which is not shared by other traditional methods of nonparametric intensity estimation such as kernel estimators and orthogonal series expansions of the intensity rather than the log-intensity. To assess the quality of the estimation, we will measure the discrepancy between an estimator $\hat{f}_t$ and the true intensity function $f$ in the sense of relative entropy (Kullback–Leibler distance) between two point processes,

$$\Delta(f; \hat{f}_t) = \int \left( f \log\left(\frac{f}{\hat{f}_t}\right) - f + \hat{f}_t \right) d\mu,$$

where the logarithms above are taken with base $e$. One can show (see Lemma VII.3 of [3]) that the above distance of the two intensities is also the Kullback–Leibler distance between the corresponding Poisson processes. It is well known that $\Delta$ is nonnegative and equals zero if and only if $\hat{f}_t = f$ a.e. The intensity $\hat{f}_t$ in the exponential family that is closest to $f$ in this relative entropy sense is the so-called information projection of $f$ [6].

The ill-posed nature of the problem comes from the assumption that $K$ is compact and therefore its inverse is not $L^2$-bounded. As in [5] we will express the ill-posed condition of $K$ by a smoothing action: $K$ will map $L^2(E, \mu)$ into some smoothness space $H^r$ for some $r > 0$. Following the notation in the paper cited above, we will say that $K$ has the smoothing property of order $\nu > 0$ if $K$ maps the Sobolev space $H^s$ onto $H^{s+\nu}$ or the Besov space $B^s_{p,q}$ onto $B^{s+\nu}_{p,q}$. Recall that the Sobolev space $H^s(\mathbb{R})$, $s \in \mathbb{R}$, is the space of tempered distributions $v$ such that

$$\|v\|_s^2 = \int_{\mathbb{R}} (1 + |\xi|^2)^s |\hat{v}(\xi)|^2 \, d\xi < \infty,$$

where

$$\hat{v}(\xi) = \int_{\mathbb{R}} e^{i\xi t} v(t) \, dt$$

denotes the Fourier transform of $v$. The Besov spaces form another particular family of smoothness spaces. Essentially the Besov spaces $B^s_{p,q}(\mathbb{R}^d)$ consist of functions that "have $s$ derivatives in $L^p$"; the parameter $q$ provides some additional fine-tuning to the definition of these spaces.

For a self-adjoint positive definite operator $K$, the smoothing property can be expressed by the ellipticity property,

$$(2.1) \qquad \langle Kf, f \rangle \sim \|f\|^2_{H^{-\nu/2}},$$

where $\langle \cdot, \cdot \rangle$ denotes the standard inner product on $L^2(E, \mu)$ and $H^{-\nu}$ is the dual space of $H^\nu$ appended with appropriate boundary conditions depending on the problem (homogeneous, periodic, etc.) (see [5]).

As already explained in the Introduction, a key ingredient for solving the Poisson inverse problem is the use of standard wavelet bases of $L^2(E, \mu)$ which allow the characterization of the function spaces that describe both the smoothness of the solution and the smoothing action of the operator $K$, since wavelet bases provide also an unconditional basis for a variety of other useful Banach spaces of functions, such as Hölder spaces, Sobolev spaces and, more generally, Besov spaces. Assume that we have a scaling function $\phi$ and a wavelet function $\psi$. Scaling and wavelet functions at scale $j$ (i.e., resolution level $2^j$) will be denoted by $\phi_\lambda$ and $\psi_\lambda$, where the index $\lambda$ summarizes both the usual scale and space parameters $j$ and $k$ [e.g., for one-dimensional wavelets, $\lambda = (j, k)$ and $\psi_{j,k} = 2^{j/2}\psi(2^j \cdot -k)$]. If $d \geq 2$, the notation $\psi_\lambda$ stands for the adaptation of scaling and wavelet functions to multidimensional domains. The notation $|\lambda| = j$ will be used to denote a wavelet at scale $j$, while $|\lambda| < j$ denotes some wavelet at scale $j'$, with $0 \leq j' < j$ (we shall assume, merely for notational convenience, that the usual coarse level of approximation $j_0$ is equal to 0). With this notation, we assume that:

(a) The scaling functions $(\phi_\lambda)_{|\lambda|=j}$ span a finite-dimensional space $V_j$ within a multiresolution hierarchy $V_0 \subset V_1 \subset \cdots \subset L^2(E, \mu)$, such that $\dim(V_j) = 2^{jd}$ (periodic wavelets for notational convenience).

(b) We are in the orthonormal case, that is, the scaling functions $(\phi_\lambda)_{|\lambda|=j}$ are an orthonormal basis of $V_j$, and the wavelets $(\psi_\lambda)_{|\lambda|=j}$ form an orthonormal basis of $W_j$ which is the orthogonal complement of $V_j$ into $V_{j+1}$.

(c) For any $g \in L^2(E, \mu)$, its wavelet decomposition can be written as

$$g = \sum_{|\lambda|=0} \tau_\lambda \phi_\lambda + \sum_{j \geq 0} \sum_{|\lambda|=j} \beta_\lambda \psi_\lambda,$$

where $\tau_\lambda = \langle g, \phi_\lambda \rangle$ and $\beta_\lambda = \langle g, \psi_\lambda \rangle$.

(d) To simplify the notation we shall use the convenient slight abuse of notation that sweeps up the coarsest-$j$ scaling functions into the $\psi_\lambda$ as well, that is, we will sometimes write $(\psi_\lambda)_{|\lambda|=-1}$ for $(\phi_\lambda)_{|\lambda|=0}$. We thus denote the complete $d$-dimensional, inhomogeneous wavelet basis by $\{\psi_\lambda; \lambda \in \Lambda\}$. By truncating the wavelet decomposition at level $j$, we obtain the orthogonal projection onto $V_j$,

$$P_j g = \sum_{|\lambda| < j} \beta_\lambda \psi_\lambda.$$

(e) We also assume that $\|\psi_\lambda\|_\infty = \|\psi\|_\infty 2^{|\lambda|d/2}$.

Wavelets provide unconditional bases for the Besov spaces, and one can express whether or not a function $g$ on $E$ belongs to a Besov space by a fairly simple and completely explicit requirement on the absolute value of the wavelet coefficients of $g$. More precisely, let us assume that the original one-dimensional $\phi$ and $\psi$ are in $C^L(\mathbb{R})$, with $L > s$, that $\sigma = s + d(1/2 - 1/p) \geq 0$, and define the norm $\| \cdot \|_{s,p,q}$ by

$$\|g\|_{s,p,q} = \left( \sum_{j=0}^{\infty} \left( 2^{j\sigma p} \sum_{\lambda \in \Lambda, |\lambda| = j} |\langle g, \psi_\lambda \rangle|^p \right)^{q/p} \right)^{1/q}.$$

Then this norm is equivalent to the traditional Besov norm, that is, there exist strictly positive constants $A$ and $B$ such that

$$A\|g\|_{s,p,q} \leq \|g\|_{B^s_{p,q}} \leq B\|g\|_{s,p,q}.$$

The condition that $\sigma \geq 0$ is imposed to ensure that $B^s_{p,q}(\mathbb{R}^d)$ is a subspace of $L^2(\mathbb{R}^d)$; we shall restrict ourselves to this case in this paper.

To end this section, and since our estimation procedures will be based on a wavelet-Galerkin projection method, we recall here some useful results on linear Galerkin projection methods for solving linear problems $h = Kf$. For a more detailed description the reader is referred to the the fairly extensive presentation in the paper by Cohen, Hoffmann and Reiss [5].

Let $f \in L^2(E, \mu)$; then the function $f_j \in V_j$ is said to be the Galerkin approximation of $f$ if for all $v \in V_j$

$$\langle K f_j, v \rangle = \langle K f, v \rangle.$$

Let $F_j \in \mathbb{R}^{2^{jd}}$ be the vector of wavelet coefficients of $f_j \in V_j$; then the Galerkin projection method for approximating $f$ amounts to solving the linear system

$$K_j G_j = G^K,$$

where $K_j = (\langle K \psi_\lambda, \psi_\kappa \rangle)_{|\lambda| < j, |\kappa| < j}$ is a symmetric positive definite matrix and $G^K = (\langle K f, \psi_\kappa \rangle)_{|\kappa| < j}$ is a "data" vector. Now, define the Galerkin wavelets $u^j_\lambda \in V_j$ as

$$(2.2) \qquad \langle K u^j_\lambda, v \rangle = \langle \psi_\lambda, v \rangle \qquad \text{for all } v \in V_j.$$

Let $U_\lambda^j$ be the vector of wavelet coefficients of $u_\lambda^j \in V_j$; then

$$U_\lambda^j = K_j^{-1} \Psi_\lambda,$$

where $\Psi_\lambda = (\langle \psi_\lambda, \psi_\kappa \rangle)_{|\kappa| < j}$ is a vector with zero entries except for the $\lambda$th component which is equal to 1. Note that

$$\begin{aligned}
\langle u_\lambda^j, Kf \rangle &= (U_\lambda^j)^T (\langle \psi_\kappa, Kf \rangle)_{|\kappa| < j} \\
&= (U_\lambda^j)^T G^K \\
&= \Psi_\lambda^T K_j^{-1} G^K = \Psi_\lambda^T G_j = G_{j,\lambda},
\end{aligned}$$

where $G_{j,\lambda} = \langle f_j, \psi_\lambda \rangle$ denotes the $\lambda$th component of $G_j$. Hence, if we define $f_j \in V_j$ by $\langle f_j, \psi_\lambda \rangle = \langle Kf, u_\lambda^j \rangle$, then $f_j$ is the *Galerkin approximation* of $f$.

**3. Information projection-based estimation.** Information projection for the estimation of density functions has been studied by Barron and Sheu [2]. They obtained various existence results and asymptotic bounds for the distance $\int p \log(p/q)$ between two probability density functions $p$ and $q$. Their estimation procedure is based on sequences of exponential families spanned by orthogonal functions such as polynomials, splines and trigonometric series. Estimation of density functions by approximation of log-densities with wavelets has been considered by Koo and Kim [17].

We adapt in this section the results of Barron and Sheu [2] to the context of log-intensity functions approximated by wavelet series with the use of the Kullback–Leibler distance between two point processes. More precisely, let $j \geq 0$. If $\theta$ denotes a vector in $\mathbb{R}^{2^{jd}}$, then $\theta_\lambda$ denotes its $\lambda$th component. The wavelet-based exponential family $\mathcal{E}_j$ at scale $j$ will be defined as the set of functions

$$\mathcal{E}_j = \left\{ f_{j,\theta}(\cdot) = \exp\left( \sum_{|\lambda| < j} \theta_\lambda \psi_\lambda(\cdot) \right), \theta = (\theta_\lambda)_{|\lambda| < j} \in \mathbb{R}^{2^{jd}} \right\}.$$

Following Csiszár [6], the intensity $f_{j,\theta}$ in the exponential family $\mathcal{E}_j$ that is closest to the true intensity $f$ in the relative entropy sense is characterized as the unique intensity function in the family for which $\langle f_{j,\hat{\theta}_t}, \psi_\lambda \rangle = \langle f, \psi_\lambda \rangle$. It seems therefore natural to estimate the unknown intensity function $f$ by searching for some $\hat{\theta}_t \in \mathbb{R}^{2^{jd}}$ such that

$$\langle f_{j,\hat{\theta}_t}, \psi_\lambda \rangle = \frac{1}{t} \int u_\lambda^j \, dG = \hat{\alpha}_\lambda^t \qquad \text{for all } |\lambda| < j.$$

If there exists a solution to this problem, then $f_{j,\hat{\theta}_t}$ will be called the *Galerkin information projection* estimate of $f$ at scale $j$, since in the context of the wavelet-Galerkin approach for solving $y = Kf + \sigma \, dW$, the estimation $\langle y, u_\lambda^j \rangle$ of $\langle Kf, u_\lambda^j \rangle$ is replaced by $\frac{1}{t} \int u_\lambda^j \, dG = \hat{\alpha}_\lambda^t$, while $\langle f_j, \psi_\lambda \rangle$ is replaced by $\langle f_{j,\hat{\theta}_t}, \psi_\lambda \rangle$.

We already pointed out the advantage of such an approach since one can guarantee that the intensity function estimates are positive. The following lemma states some of the I-projection properties onto $\mathcal{E}_j$ (see also [6]).

LEMMA 3.1. *Let $\alpha \in \mathbb{R}^{2^{jd}}$. Assume that there exists some $\theta(\alpha) \in \mathbb{R}^{2^{jd}}$ such that for all $|\lambda| < j$*

$$\langle f_{j,\theta(\alpha)}, \psi_\lambda \rangle = \alpha_\lambda.$$

*Then, for any intensity function $f \in L^2(E, \mu)$ such that $\langle f, \psi_\lambda \rangle = \alpha_\lambda$ and for all $\theta \in \mathbb{R}^{2^{jd}}$, the following Pythagorean-like identity holds:*

$$\Delta(f; f_{j,\theta}) = \Delta(f; f_{j,\theta(\alpha)}) + \Delta(f_{j,\theta(\alpha)}, f_{j,\theta}).$$

A consequence of the above lemma, and since $\Delta(f; h) > 0$ unless $f = h$ almost everywhere, is that $\theta(\alpha)$ (if it exists) uniquely minimizes $\Delta(f; f_{j,\theta})$ for $\theta \in \mathbb{R}^{2^{jd}}$.

From now on assume that there exists some constant $A_j < \infty$ such that for all $v \in V_j$

$$\|v\|_\infty \leq A_j \|v\|_{L^2}.$$

A key lemma relating distances between the intensities in the parametric family to distance between the corresponding wavelet coefficients is then the following.

LEMMA 3.2. *Let $\theta_0 \in \mathbb{R}^{2^{jd}}$, $\alpha_{0,\lambda} = \langle f_{j,\theta_0}, \psi_\lambda \rangle$ and $\alpha \in \mathbb{R}^{2^{jd}}$ be a given vector. Let $b = \exp(\|\log(f_{j,\theta_0})\|_\infty)$ and $e = \exp(1)$. If $\|\alpha - \alpha_0\|_2 \leq \frac{1}{2ebA_j}$, then the solution $\theta(\alpha)$ to*

$$\langle f_{j,\theta(\alpha)}, \psi_\lambda \rangle = \alpha_\lambda \qquad \text{for all } |\lambda| < j$$

*exists and satisfies*

$$(3.1) \qquad \|\theta(\alpha) - \theta_0\|_2 \leq 2eb\|\alpha - \alpha_0\|_2,$$

$$(3.2) \qquad \left\| \log\left( \frac{f_{j,\theta(\alpha_0)}}{f_{j,\theta(\alpha)}} \right) \right\|_\infty \leq 2ebA_j \|\alpha - \alpha_0\|_2,$$

$$(3.3) \qquad \Delta(f_{j,\theta(\alpha_0)}; f_{j,\theta(\alpha)}) \leq 2eb\|\alpha - \alpha_0\|_2^2.$$

The proof of this lemma relies upon a series of lemmas on bounds within exponential families for the Kullback–Leibler distance and is given in the Appendix.

**4. Linear estimation.** Let $M$ be some fixed constant and let $F_{p,q}^s(M)$ denote the set of scaled intensity functions such that

$$F_{p,q}^s(M) = \{f = \exp(g), \ \|g\|_{B_{p,q}^s} \leq M\}.$$

Note that assuming that $f \in F_{p,q}^s(M)$ implies that $f$ is strictly positive.

For $f \in F_{p,q}^s(M)$, let $g = \log_e(f)$ and define

$$D_j = \|g - P_j g\|_{L^2},$$

$$\gamma_j = \|g - P_j g\|_\infty.$$

Basic to our analysis is a decomposition of the relative entropy between the true and the estimated intensities into the sum of two terms which correspond to approximation error and estimation error (bias and variance in a familiar mean squared error analysis). The proof of this result, which relies upon some concentration inequalities for Poisson processes, is postponed to the [Appendix](#).

THEOREM 4.1.   *Assume that $\psi$ is compactly supported and that $f \in F_{p,q}^s(M)$ (with $s > d/p \geq d/2$). Let $M_1 > 1$ be a constant such that $M_1^{-1} \leq f \leq M_1$ (see Lemma* [A.4](#)*), and let $\varepsilon_j = 2M_1^2 e^{2\gamma_j+1} D_j A_j$. If $\varepsilon_j \leq 1$, the information projection exists, that is, there exists $\theta_j^* \in \mathbb{R}^{2^{jd}}$ such that*

$$\langle f_{j,\theta_j^*}, \psi_\lambda \rangle = \langle f, \psi_\lambda \rangle \qquad \text{for all } |\lambda| < j,$$

*and the approximation error satisfies*

$$\Delta(f; f_{j,\theta_j^*}) \leq C e^{\gamma_j} D_j^2.$$

*Moreover, suppose that $\psi$ is in $H^{s+d/2-\nu}$ (with $s > \nu - d/2$) and has $r$ vanishing moments with $r > s + d/2$. Let $\delta_j^t = 4M_1^2 e^{2\varepsilon_j+2\gamma_j+2} A_j^2 \rho_{j,t}$, where $\rho_{j,t} = (2^{j(\nu+d/2)}/\sqrt{t} + 2^{j(\nu+(3/2)d)}/t)^2 + 2^{-2js}$. If $\delta_j^t \leq 1$, then for every $\eta^2 \leq \frac{1}{\delta_j^t}$ there is a set of probability less than $\exp(-\eta)$, such that outside this set there exists some $\hat{\theta}_t \in \mathbb{R}^{2^{jd}}$ which satisfies*

$$\langle f_{j,\hat{\theta}_t}, \psi_\lambda \rangle = \frac{1}{t} \int u_\lambda^j \, dG \qquad \text{for all } |\lambda| < j,$$

*and the estimation error satisfies*

$$\Delta(f_{j,\theta_j^*}; f_{j,\hat{\theta}_t}) \leq C \eta^2 e^{1+\gamma_j+\varepsilon_j} \rho_{j,t}.$$

Note that by using the above theorem, explicit bounds are obtained which are applicable for each finite value of $j$ and $t$, subject to $\varepsilon_j$ and $\delta_j^t \leq 1$. We can now state the general result on the nonadaptive Galerkin information projection estimator of the unknown intensity function.

THEOREM 4.2. *Assume that $\psi$ is compactly supported and that $f \in F_{2,2}^s(M)$ (with $s > d/2$). Moreover, suppose that $\psi$ is in $H^{s+d/2-\nu}$ (with $s > \nu - d/2$) and has $r$ vanishing moments with $r > s + d/2$. Let $j(t)$ be such that $2^{-j(t)} = (\frac{1}{t})^{1/(2s+2\nu+d)}$. Then, with probability tending to 1 as $t \to \infty$, the Galerkin information projection exists and satisfies*

$$\Delta(f; f_{j(t),\hat{\theta}_t}) \leq O\left(\left(\frac{1}{t}\right)^{2s/(2s+2\nu+d)}\right).$$

The above estimator therefore converges almost surely with the optimal rate for intensities in $F_{2,2}^s(M)$. However, the main defect of the estimator defined in Theorem 4.2 is that it is suited for smooth functions and does not attain the optimal rates when, for example, $g = \log(f)$ has singularities. We therefore propose in the next section another estimator derived by applying an appropriate nonlinear thresholding procedure.

**5. Nonlinear estimation.** It is well known that linear estimators do not achieve the optimal rates of convergence when the functions to be recovered belong to Besov spaces $B_{p,p}^s$ with index $1 \leq p < 2$ (the case of functions which are not very smooth). In order to attain such a rate we need therefore some kind of nonlinear procedure and this is our aim in this section.

Our estimation procedure simply consists of applying a soft thresholding algorithm on the "data" to which we apply the Galerkin information projection inversion which was described previously, exploiting the fact that the model with Poisson intensity is not too different from the usual Gaussian white-noise model.

Let us first recall that the coefficients defining the *Galerkin information projection* estimate of $f$ at scale $j$, as derived in the previous section, are given by

$$\hat{\alpha}_{t,\lambda} = \frac{1}{t}\int u_\lambda^j\, dG = (U_\lambda^j)^T \frac{1}{t}\left(\int \psi_\mu\, dG\right)_{|\mu|<j},$$

where $U_\lambda^j = K_j^{-1}\Psi_\lambda$.

For some $j \geq 0$ (to be fixed further), we define the thresholded coefficients

$$P_t(\hat{\alpha}_{t,\lambda}) = (U_\lambda^j)^T \left(T_{\varepsilon(t)}\left(\frac{1}{t}\int \psi_\kappa\, dG\right)\right)_{|\kappa|<j} \qquad \text{for all } |\lambda| < j,$$

where $T_{\varepsilon(t)}(x) = \text{sign}(x)(x - \varepsilon(t))_+$ for $x \in \mathbb{R}$ denotes the usual soft thresholding operator with threshold $\varepsilon(t)$.

In order to build an optimal solution for the Poisson inverse problem we will use a level-dependent wavelet thresholding procedure by setting $\varepsilon(t) = t^{-1/2} 2^{\nu|\lambda|}\sqrt{|\log t|}$. The role of $2^{\nu|\lambda|}$ is to take into account the amplification of the noise by the inversion process. The following theorem shows that the resulting estimator behaves in an optimal way provided that the cutoff resolution level $j(t)$ is chosen such that $2^{-j(t)} \leq t^{-1/(2\nu)}$, where $\nu$ is the degree of ill-posedness of the estimator.

THEOREM 5.1. *Assume that $\psi$ is compactly supported and that $f \in F_{p,p}^s(M)$ with $s > 0$ and $1/p = 1/2 + s/(2\nu + d)$. Moreover, suppose that $\psi$ is in $H^{s+d/2-\nu}$ (with $s > \nu - d/2$) and has $r$ vanishing moments with $r > s + d/2$. Also assume that $K$ is an isomorphism between $L^2$ and $H^\nu$ and that it has the smoothing property of order $\nu$ with respect to the space $B_{p,p}^s$. Then, the above described Galerkin information projection estimator, say $f_{j(t),\hat{\theta}_t}$, satisfies the minimax rate*

$$\mathbb{E}(\Delta(f; f_{j(t),\hat{\theta}_t})) \leq O\left(\left(\frac{1}{t}\sqrt{|\log t|}\right)^{2s/(2s+2\nu+d)}\right),$$

*provided that $j(t)$ is such that $2^{-j(t)} \leq t^{-1/(2\nu)}$.*

Note that the lower bound on $j(t)$ does not depend on the unknown smoothness of $f$ and therefore Theorem 5.1 allows us to build an adaptive solution to our Poisson inverse problem. The assumption that $K^{-1}$ maps $H^\nu$ into $L^2$ in the above theorem is also implicit in the vaguelette–wavelet method of Donoho [7] for white-noise inverse problems.

**6. Implementation and some numerical results.** The purpose of this section is to describe the implementation of our approach and to briefly explore the performance of our method from a numerical point of view. As in [5] we will focus on a simple example of a logarithmic potential kernel in dimension 1. We will consider its action on two typical test intensity functions, which, together with their folded versions by the action of $K$, are displayed in Figure 1.

The logarithmic potential operator $K$ that we will consider is defined by

$$Kf(x) = \int_0^1 k(x,y)f(y)\,dy,$$

where

$$k(x,y) = -\log\left(\frac{1}{2}\left|\sin\frac{y-x}{2}\right|\right), \qquad x,y \in [0,1].$$

Such a kernel is singular on the diagonal $x = y$ but integrable. The corresponding operator is known to be an elliptic operator of order $-1$, which maps $H^{-1/2}$ into $H^{1/2}$ and therefore satisfies the assumptions made in this paper with $\nu = 1$. The first test function we will consider is

$$f(x) = \max\{1 - |30(x-0.5)|, 0.1\}, \qquad x \in [0,1],$$

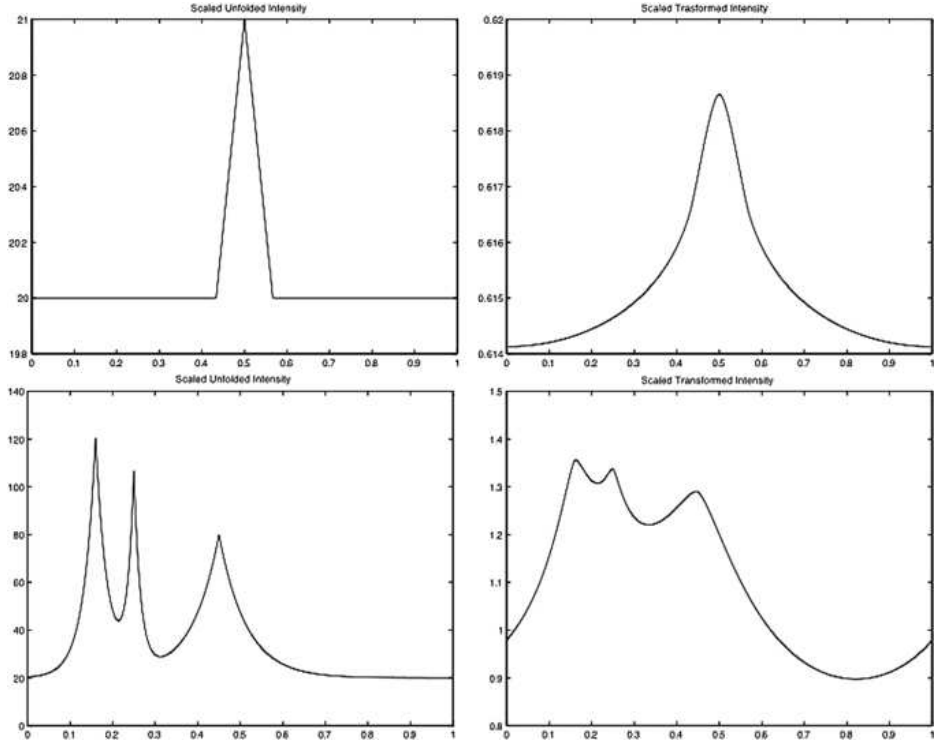

FIG. 1. *Artificial intensity functions (left) and their folded versions (right) by the action of the logarithmic potential kernel. Top-left: the intense peak; bottom-left: the burst-like intensity.*

which presents an intense peak and is badly approximated by the singular functions of $K$ but has a very sparse representation in a wavelet basis. We will also consider a fast rise-exponential decay model, giving rise to the abbreviation "FRED" in the astronomy literature, to model burst phenomena which are of the form $f(x) = f_0 + \sum_{i=1}^{3} f_i(x)$ where $f_0$ models the relatively constant background level of gamma-ray photon arrivals and $f_i$ models the $i$th peak in the burst of a form

$$f_i(x) = \begin{cases} a_i \exp(-|x - m_i|/\sigma_{r,i}^{\nu_i}), & \text{if } x \leq m_i, \\ a_i \exp(-|x - m_i|/\sigma_{d,i}^{\nu_i}), & \text{if } x > m_i. \end{cases}$$

In the expression above $m_i$ denotes the location of the $i$th peak, and the factors $a_i$, $\sigma_{r,i}$, $\sigma_{d,i}$ and $\nu_i$ control respectively the amplitude, the rise, the decay and the peakedness.

Most often data can only be observed in binned form because of the discrete nature of the measurement apparatus or because binning may be enforced by data handling and computing efficiency. We therefore have cho-sen to discretize $K$ for a maximal resolution level $J = 11$ by computing the stiffness matrix $K_J$ with entries

$$(K_J)_{\ell,k=0,\dots,2^J-1} = (\langle K_J\phi_{J,\ell}, \phi_{J,k}\rangle)_{\ell,k=0,\dots,2^J-1},$$

where the $\phi_{J,k} = 2^{J/2}I_{[k2^{-J},(k+1)2^{-J}]}$ are the Haar scaling functions. Each integral

$$\langle K_J\phi_{J,\ell}, \phi_{J,k}\rangle = \int_0^1 \int_0^1 k(x,y)\phi_{J,\ell}(x)\phi_{J,k}(y)\,dx\,dy$$

was computed by Riemann approximation at scale $2^{-16}$. Note that the kernel $k$ is such that $k(x,y) = h(x-y)$ where $h(\cdot)$ is a 1-periodic function. The discretization $K_J$ of $K$ is therefore a Toeplitz cyclic matrix and the fast Fourier transform makes the computation of the action of $K$ on functions approximated in the Haar basis numerically fast and easy. But this is not the only reason we have used the Haar-based discretization for both $f$ and $K_J$. More specifically, the Haar scaling basis at resolution $J$ induces a partition of the interval $[0,1)$ into $2^J$ disjoint and measurable bins $B_k = [k2^{-J},(k+1)2^{-J})$. Integrating the function $tK_Jf$ with respect to the Poisson counting measure $G$ simply leads to observed data consisting of counts observed in the bins $B_k$. By the Poisson nature of $G$, these are independent Poisson counts with expected values within each bin $t\int_{B_k} h$, $k = 0,\dots,2^J-1$. Moreover, taking a high-resolution $J$ permits the approximation of $\int_{B_k} h$ by $2^{-J}h(k/2^J)$ and this is what we have done for creating the simulated data in the examples.

For the examples treated in this paper, the estimation was implemented using Symmlets with six vanishing moments. Since the set $\{x_1,\dots,x_n\}$ of the $n = 2^J$ points at which the data is sampled is dyadic, any scalar product involving a wavelet at a lower resolution is computed via the discrete wavelet transform. The information projection estimator was obtained by solving the system of equations given in Theorem 5.1 at some maximal resolution $J_{\max}$.

To find the estimate $\hat{\boldsymbol{\theta}}_t$ we have used, inspired by a similar approach in [13], a modified version of the Newton–Raphson method. Let $S(\boldsymbol{\theta})$ denote the $J_{\max}$-dimensional vector of elements

$$S(\boldsymbol{\theta})_\lambda = (P_t(\hat{\alpha}_{t,\lambda}) - \langle f_{t,\boldsymbol{\theta}}, \psi_\lambda\rangle),$$

and $H(\boldsymbol{\theta})$ the $J_{\max} \times J_{\max}$ Hessian matrix whose $(\lambda, \lambda')$ entry is given by

$$H(\boldsymbol{\theta})_{\lambda,\lambda'} = \int f_{t,\boldsymbol{\theta}}(x)\psi_\lambda(x)\psi_{\lambda'}(x)\,dx.$$

The method to compute $\hat{\boldsymbol{\theta}}_t$ is to start with an initial guess $\boldsymbol{\theta}^0$ and iteratively determine $\boldsymbol{\theta}^{m+1}$ according to

$$\boldsymbol{\theta}^{m+1} = \boldsymbol{\theta}^m + H^{-1}(\boldsymbol{\theta}^m)S(\boldsymbol{\theta}^m),$$

with a standard criterion for stopping the iterations.

One difficulty in implementing the above algorithm is that Symmlets have no closed-form functional expression and the integration involved in the computation of the Hessian $H$ can be very time consuming if one has to compute a table of the appropriate values of the function $\psi_\lambda$. We have used instead an efficient filter-bank algorithm for computing such integrals similar to the one used by Vannucci and Corradi [27] or Kovac and Silverman [18] to compute the diagonal elements of the covariance structure of the wavelet coefficients, which amounts to computing the fast two-dimensional wavelet transform of the diagonal matrix whose diagonal is the vector $f_{t,\boldsymbol{\theta}}(x_i)$, $i = 1, \dots, J_{\max}$, and to retain only the diagonal blocks of the transform.

For our first example, we consider the peaky function and choose a maximal resolution $J_{\max} = 10$ and an "observation time" $t = 10^8$ (corresponding to the noise level used by Cohen, Hoffmann and Reiss [5] for a similar white-noise model). A typical sample from the simulated model is shown in the top-left panel of Figure 2.

Let us recall that our nonlinear information projection estimator depends on the cutoff level $j(t)$ given by $2^{-j(t)} \le (\frac{1}{t})^{1/(2\nu)}$ and the level-dependent thresholds $\varepsilon(t) = 2^{\nu|\lambda|}t^{-1/2}\sqrt{|\log t|}$. We therefore have used these expressions with $\nu = 1$ to estimate the unfolded intensity function. The top-right panel of Figure 2 displays the nonlinear Galerkin estimator for the folded Poisson data displayed in the top-left panel of Figure 2. We observe, and this is true also for the second example, that the peak is very well estimated. However, some oscillations are observed on the right side of the central peak. A possible remedy to this defect could be to use a translation-invariant procedure, but such an approach is beyond the scope of this paper.

Our second example concerns the burst-like intensity function. The data displayed in the bottom-left panel of Figure 2 is simulated as above using for an unfolded intensity a burst signal with a constant intensity of 20 assigned to the background and three peaks. Since we are using the same logarithmic potential kernel to fold the intensity, we have also used here a value of $\nu$ equal to 1. Maximal resolution, smooth cutoff level and thresholds were chosen exactly as in the previous example and the estimation procedure provides us the fit in the bottom-right panel of Figure 2 for the unfolded burst-like intensity, confirming the good behavior of our procedure even for complicated intensity structures.

**7. Conclusions.** The methodology of this paper was motivated by wavelet–vaguelette decomposition (WVD) methods that have been developed in the literature for solving inverse problems with Gaussian white-noise perturbations. Such methods are most appropriate only for the restricted class of homogeneous operators, particularly from a computational perspective, and extra theoretical and numerical work is required to handle more general operators or Poisson inverse problems. The method developed is particularly well suited for our Poisson problem because it is designed for positive definite operators, its numerical implementation is straightforward, it can easily be extended to approximately known operators and it leads to positive estimated intensities. It combines the numerical simplicity of Galerkin projection methods on a high-dimensional space as an inversion procedure and wavelet thresholding as an adaptive smoothing technique.

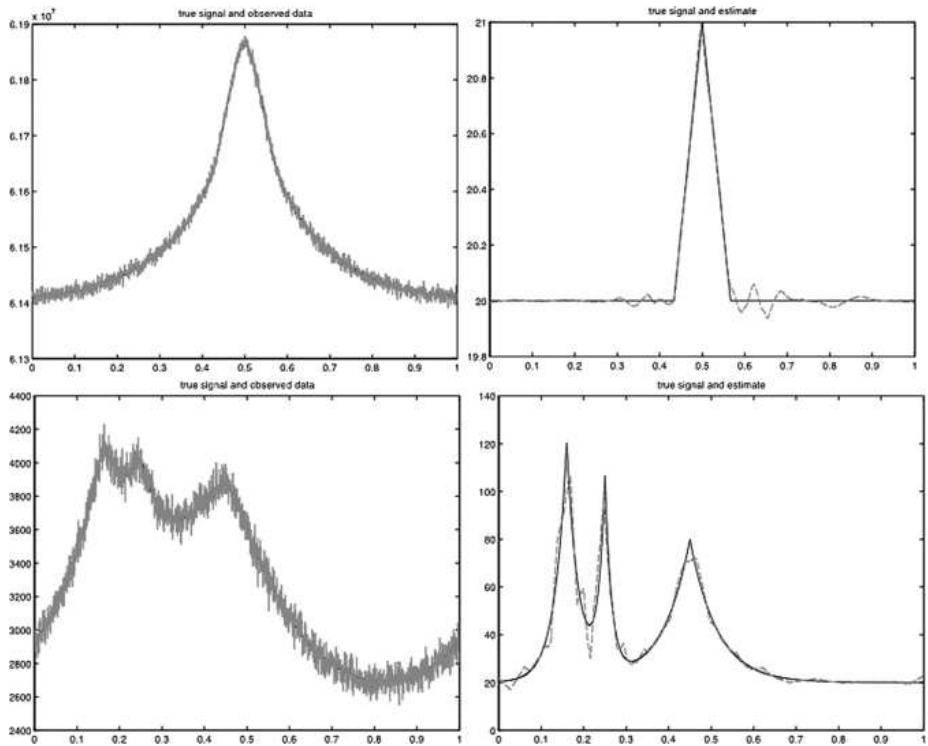

FIG. 2. *Simulated Poisson data obtained by folding the intensity functions with the logarithmic potential. Top-left: the intense peak; bottom-left: the burst-like intensity. The right panels display the corresponding unknown intensities (solid) and their nonlinear Garlekin estimates (dashed).*

Our attention was restricted to inverse problems in which $E_0$ and $E_1$ are identical Borel subsets of $\mathbb{R}^d$ and the operator $K$ is self-adjoint positive definite. In the case where $E_0 \neq E_1$ and $K$ is not self-adjoint but just injective, one may choose, as is done in [5], $E_1 = K(E_0)$ and solve instead the inverse problem $K^* h = K^* K f$ where $K^*$ denotes the adjoint of $K$. In such a way we are led back to the wavelet-Galerkin information projection method developed in this paper.

Another approach, which is worth investigating in the future, would be to derive a Landweber-type iterative algorithm that involves a denoising procedure at each iteration step and provides a sequence of approximations converging in norm to the maximum penalized likelihood minimizer as is done by Figueiredo and Nowak [10], who derive such an iterative algorithm for inverting a convolution operator acting on objects that are sparse in the wavelet domain.

## APPENDIX: PROOFS OF THE MAIN RESULTS

We first derive the proof of Lemma 3.1.

PROOF OF LEMMA 3.1. Note that

$$\Delta(f; f_{j,\theta}) = D(f, f_{j,\theta}) + \int (-f + f_{j,\theta(\alpha)}) \, d\mu + \int (-f_{j,\theta(\alpha)} + f_{j,\theta}) \, d\mu,$$

where

$$D(f, f_{j,\theta}) = \int f \log\left(\frac{f}{f_{j,\theta}}\right) d\mu = D(f, f_{j,\theta(\alpha)}) + D(f_{j,\theta(\alpha)}, f_{j,\theta})$$

by Lemma 4 of [2], which completes the proof.  □

To prove Lemma 3.2 we need some preliminary lemmas on bounds within exponential families for the Kullback–Leibler distance; their proofs can be easily derived from the proofs of Lemmas 1 and 4 of [2] and thus they are omitted.

LEMMA A.1.  *Let $f$ and $h$ be two intensity measures in $L^2(E, \mu)$. Assume that $\log(\frac{f}{h})$ is bounded; then*

$$\Delta(f; h) \geq \frac{1}{2} e^{-\|\log(f/h)\|_\infty} \int f\left(\log\left(\frac{f}{h}\right)\right)^2 d\mu,$$

$$\Delta(f; h) \leq \frac{1}{2} e^{\|\log(f/h)\|_\infty} \int f\left(\log\left(\frac{f}{h}\right)\right)^2 d\mu.$$

LEMMA A.2.  *For $j \geq 0$, let $\theta_0, \theta \in \mathbb{R}^{2dj}$ and $b = \exp(\|\log(f_{j,\theta_0})\|_\infty)$; then*

$$\left\|\log\left(\frac{f_{j,\theta_0}}{f_{j,\theta}}\right)\right\|_\infty \leq A_j \|\theta_0 - \theta\|_2,$$

$$\Delta(f_{j,\theta_0}; f_{j,\theta}) \geq \frac{1}{2b} e^{-A_j\|\theta_0 - \theta\|_\infty} \|\theta_0 - \theta\|_2^2,$$

$$\Delta(f_{j,\theta_0}; f_{j,\theta}) \leq \frac{b}{2} e^{A_j\|\theta_0 - \theta\|_\infty} \|\theta_0 - \theta\|_2^2.$$

We can proceed to the proof of Lemma 3.2.

PROOF OF LEMMA 3.2.    This proof is inspired by the proof of Lemma 5 of [2]. Let

$$F(\theta) = \theta \cdot \alpha - H(\theta),$$

where $H(\theta) = \int f_{j,\theta}(x) \, d\mu(x)$. If $\alpha = \alpha_0$, then the result is trivial. Now, if $\alpha \neq \alpha_0$, note that for any $\theta \in \mathbb{R}^{2dj}$

$$\Delta(f_{j,\theta_0}; f_{j,\theta}) = \alpha_0 \cdot (\theta_0 - \theta) + H(\theta) - H(\theta_0).$$

Hence,

$$F(\theta_0) - F(\theta) = \Delta(f_{j,\theta_0}; f_{j,\theta}) - (\alpha_0 - \alpha) \cdot (\theta_0 - \theta).$$

So, by Lemma A.2 and the Cauchy–Schwarz inequality, we have that

$$F(\theta_0) - F(\theta) \geq \frac{1}{2b} e^{-A_j \|\theta_0 - \theta\|_\infty} \|\theta_0 - \theta\|_2^2 - \|\alpha_0 - \alpha\|_2 \|\theta_0 - \theta\|_2.$$

This inequality is strict if $\theta \neq \theta_0$. For all $\theta$ such that $\|\theta_0 - \theta\|_2 = 2eb\|\alpha_0 - \alpha\|_2$,

$$F(\theta_0) - F(\theta) > 2eb\|\alpha_0 - \alpha\|_2^2 (e^{1-2A_j eb\|\alpha_0 - \alpha\|_2} - 1).$$

The right-hand side is positive whenever $2A_j eb\|\alpha_0 - \alpha\|_2 \leq 1$. Hence, $F(\theta_0) > F(\theta)$ for all $\theta$ such that $\|\theta_0 - \theta\|_2 = 2eb\|\alpha_0 - \alpha\|_2$. Consequently, $F$ has an extreme point $\theta^*$ such that $\|\theta_0 - \theta^*\|_2 < 2eb\|\alpha_0 - \alpha\|_2$. The gradient of $F$ at $\theta^*$ must satisfy

$$\langle f_{j,\theta^*}, \psi_\lambda \rangle = \alpha_\lambda \qquad \text{for all } |\lambda| < j,$$

and so $\theta(\alpha) = \theta^*$. Hence, inequality (3.1) immediately follows. Inequality (3.2) follows from Lemma A.2. Since $F(\theta(\alpha)) \geq F(\theta_0)$, we have that

$$\begin{aligned}
\Delta(f_{j,\theta(\alpha_0)}; f_{j,\theta(\alpha)}) &\leq (\alpha_0 - \alpha) \cdot (\theta_0 - \theta(\alpha)) \\
&\leq \|\alpha_0 - \alpha\|_2 \|\theta_0 - \theta\|_2 \\
&\leq 2eb\|\alpha_0 - \alpha\|_2^2,
\end{aligned}$$

which completes the proof.  □

To prove the main results of this paper we shall need a series of technical lemmas, stated and proved below. Throughout this section, $C$ will denote a constant whose value may change from line to line.

LEMMA A.3.   *If $\psi$ is compactly supported, then*

$$\left\| \sum_{|\lambda|=j} \beta_\lambda \psi_\lambda \right\|_\infty \leq C 2^{jd/2} \|\beta_j\|_2.$$

PROOF. This lemma immediately follows from the proof of Lemma 1 of [17]. □

The following lemma is similar to Lemma 2 of [17].

LEMMA A.4. *Assume that $\psi$ is compactly supported and that $f \in F_{p,q}^s(M)$. If $s > d/p \geq d/2$, then there exists a constant $M_1$ such that*

$$0 < \frac{1}{M_1} \leq f \leq M_1 < \infty.$$

PROOF. Let $g = \log(f) = \sum_{j=-1}^{\infty} \sum_{|\lambda|=j} \beta_\lambda \psi_\lambda$. Since $\|g\|_{B_{p,q}^s} \leq M$, we have

$$\|\beta_j\|_p = \left( \sum_{|\lambda|=j} |\beta_\lambda|^p \right)^{1/p} \leq M 2^{-js'},$$

where $s' = s + d(1/2 - 1/p)$. If $s > d/p \geq d/2$, then

(A.1) $$\|\beta_j\|_2 \leq \|\beta_j\|_p \leq C 2^{-js'}.$$

Then, by Lemma A.3

$$\begin{aligned}
\|g\|_\infty &\leq \sum_{j=-1}^{\infty} \left\| \sum_{|\lambda|=j} \beta_\lambda \psi_\lambda \right\|_\infty \\
&\leq \sum_{j=0}^{\infty} C 2^{jd/2} \|\beta_j\|_2 \\
&\leq C \sum_{j=0}^{\infty} 2^{j(d/2 - s')} \\
&\leq C \sum_{j=0}^{\infty} 2^{-j(s-d/p)}.
\end{aligned}$$

Since $s > d/p$, $\sum_{j=0}^{\infty} 2^{-j(s-d/p)} < \infty$ and therefore there exists some constant $M_1 > 1$ such that $\|g\|_\infty = \|\log f\|_\infty \leq \log M_1$. □

Now we give bounds for $A_j$, $D_j$ and $\gamma_j$.

LEMMA A.5. *Assume that $\psi$ is compactly supported; then*

$$A_j = C 2^{jd/2}.$$

*Moreover suppose that $f \in F_{p,q}^s(M)$. If $s > d/p \geq d/2$, then*

$$D_j \leq C 2^{-j(s+d(1/2-1/p))},$$

$$\gamma_j \leq C 2^{-j(s-d/p)}.$$

PROOF.   The result for $A_j$ immediately follows from Lemma A.3. Note that from (A.1),

$$D_j^2 = \sum_{j' \geq j} \|\beta_{j'}\|_2^2 \leq C \sum_{j' \geq j} 2^{-2j'(s+d(1/2-1/p))} = O(2^{-2j(s+d(1/2-1/p))}).$$

By definition,

$$\gamma_j = \|g - P_j g\|_\infty \leq A_j D_j \leq C 2^{-j(s-d/p)},$$

which completes the proof.   □

We may now proceed to the proofs of our main results.

PROOF OF THEOREM 4.1.

*Approximation error term.*   Let $g = \log(f) = \sum_{j=-1}^\infty \sum_{|\lambda|=j} \beta_\lambda \psi_\lambda$, and for all $|\lambda| < j$, let $\alpha_{j,\lambda} = \langle \exp(P_j g), \psi_\lambda \rangle$ and $\alpha_\lambda = \langle f, \psi_\lambda \rangle$. Observe that the coefficients $(\alpha_{j,\lambda} - \alpha_\lambda)$, $|\lambda| < j$, are the coefficients of the orthogonal projection of $f - \exp(P_j g)$ onto $V_j$. Hence by Bessel's inequality

$$\|\alpha_j - \alpha\|_2^2 \leq \|f - \exp(P_j g)\|_{L^2}^2.$$

Given our assumptions on $\psi$ and $f$, Lemma A.4 implies that

$$\|\alpha_j - \alpha\|_2^2 \leq M_1 \int \frac{(f - \exp(P_j g))^2}{f} \, d\mu,$$

and so by Lemma 2 of [2],

$$\|\alpha_j - \alpha\|_2^2 \leq M_1 e^{2\|g - P_j g\|_\infty} \int f(g - P_j g)^2 \, d\mu \leq M_1^2 e^{2\gamma_j} D_j^2.$$

Now, apply Lemma 3.2 with $\theta_{0,\lambda} = \beta_\lambda$, $\alpha_\lambda = \langle f, \psi_\lambda \rangle$ for all $|\lambda| < j$ and $b = e^{\|\log(\exp(P_j g))\|_\infty}$.   Since   $\|\log(f/\exp(P_j g))\|_\infty = \gamma_j$,   we   have   that $\|\log(\exp(P_j g))\|_\infty \leq \log M_1 + \gamma_j$ and therefore $b \leq M_1 e^{\gamma_j}$. From Lemma 3.2, we have that if $M_1 e^{\gamma_j} D_j \leq \frac{1}{2ebA_j}$, that is, if $\varepsilon_j \leq 1$, then $\theta_j^* = \theta(\alpha)$ exists. By Lemma 3.1 (Pythagorean-like relationship), we obtain that

$$\Delta(f; f_{j,\theta_j^*}) \leq \Delta(f; \exp(P_j g)).$$

Thence, by Lemma A.1,

$$\Delta(f; f_{j,\theta_j^*}) \leq \tfrac{1}{2} e^{\|g - P_j g\|_\infty} \int f(g - P_j g)^2 \, d\mu$$
$$\leq \tfrac{1}{2} M_1 e^{\gamma_j} D_j^2,$$

which completes the proof of the first assertion of the theorem.

*Estimation error term.* Using the above notation and proof, and since by assumption $\varepsilon_j \leq 1$, let $\theta_j^* \in \mathbb{R}^{2jd}$ be the parameter vector achieving the minimum of the relative entropy for intensities in the exponential family. For all $|\lambda| < j$, define $\alpha_{0,\lambda} = \langle f, \psi_\lambda \rangle = \langle f_{j,\theta_j^*}, \psi_\lambda \rangle$ and let $\hat{\alpha}_{t,\lambda} = \frac{1}{t} \int u_\lambda^j \, dG$. It is easy to see that $\mathbb{E}(\hat{\alpha}_{t,\lambda}) = \langle u_\lambda^j, Kf \rangle$. We now have

$$
\begin{aligned}
\|\hat{\alpha}_t - \alpha_0\|_2^2 &= \sum_{|\lambda| < j} (\hat{\alpha}_{t,\lambda} - \alpha_{0,\lambda})^2 \\
&\leq 2 \left( \sum_{|\lambda| < j} (\hat{\alpha}_{t,\lambda} - \langle u_\lambda^j, Kf \rangle)^2 + (\langle u_\lambda^j, Kf \rangle - \alpha_{0,\lambda})^2 \right).
\end{aligned}
\tag{A.2}
$$

Concerning the second term of the right-hand side of inequality (A.2), using expression (2.2), note that

$$
\begin{aligned}
(\langle u_\lambda^j, Kf \rangle - \langle f, \psi_\lambda \rangle)^2 &= (\langle f, Ku_\lambda^j - \psi_\lambda \rangle)^2 \leq \|f\|_{L^2}^2 \|Ku_\lambda^j - \psi_\lambda\|_{L^2}^2 \\
&= \|f\|_{L^2}^2 (\|Ku_\lambda^j\|_{L^2}^2 + 1 - 2\langle Ku_\lambda^j, \psi_\lambda \rangle).
\end{aligned}
$$

By definition we have that $\langle Ku_\lambda^j, \psi_\lambda \rangle = \langle \psi_\lambda, \psi_\lambda \rangle = 1$. It follows that

$$
\begin{aligned}
\|Ku_\lambda^j\|_{L^2}^2 &= \sum_{|\mu| < j} \langle Ku_\lambda^j, \psi_\mu \rangle^2 + \sum_{|\mu| \geq j} \langle Ku_\lambda^j, \psi_\mu \rangle^2 \\
&= 1 + \sum_{|\mu| \geq j} \langle Ku_\lambda^j, \psi_\mu \rangle^2.
\end{aligned}
$$

Given the assumptions on the wavelet $\psi$ and the operator $K$, $u_\lambda^j$ belongs to $H^{s+d/2-\nu}$ and hence $Ku_\lambda^j$ belongs to $H^{s+d/2}$. Since $r > s + d/2$, it follows from approximation theory that

$$
\sum_{|\mu| \geq j} \langle Ku_\lambda^j, \psi_\mu \rangle^2 \leq 2^{-2j(s+d/2)},
$$

and therefore

$$
(\langle u_\lambda^j, Kf \rangle - \langle f, \psi_\lambda \rangle)^2 \leq \|f\|_{L^2}^2 2^{-2js - jd}.
$$

Thence we obtain for the second term of the right-hand side of inequality (A.2)

$$
\sum_{|\lambda| < j} (\langle u_\lambda^j, Kf \rangle - \langle f, \psi_\lambda \rangle)^2 \leq \|f\|_{L^2}^2 2^{-2js}.
\tag{A.3}
$$

To control the first term of the right-hand side of (A.2), let $S_j = \mathrm{span}\{u_\lambda^j; |\lambda| < j\}$ and set

$$
\chi^2(S_j) = \sum_{|\lambda| < j} (\hat{\alpha}_{t,\lambda} - \langle u_\lambda^j, Kf \rangle)^2 = \sum_{|\lambda| < j} \left( \frac{1}{t} \int u_\lambda^j \, dG - \langle u_\lambda^j, Kf \rangle \right)^2.
$$

Noticing that

$$\chi(S_j) = \sup_{\{a \in \mathbb{R}^{jd}; \|(a_\lambda)_{|\lambda|<j}\|_2 < 1\}} \int \frac{\sum_{|\lambda|<j} a_\lambda u_\lambda^j}{t}(dG - tKf\,d\mu),$$

we can use Corollary 2 of [24] about concentration inequalities for Poisson processes to get, for any $w > 0$ and $\varepsilon > 0$,

$$(A.4) \quad \mathbb{P}\{\chi(S_j) \geq (1+\varepsilon)\mathbb{E}(\chi(S_j)) + \sqrt{12v_0 w} + \kappa(\varepsilon)b_0 w\} \leq \exp(-w),$$

where $\kappa(\varepsilon) = 1.25 + 32/\varepsilon$,

$$v_0 = \sup_{\{a \in \mathbb{R}^{jd}; \|a\|_2 < 1\}} \int \frac{\sum_{|\lambda|<j}(a_\lambda u_\lambda^j)^2}{t^2} tKf\,d\mu$$

and

$$b_0 = \sup_{\{a \in \mathbb{R}^{jd}; \|a\|_2 < 1\}} \frac{\|\sum_{|\lambda|<j} a_\lambda u_\lambda^j\|_\infty}{t}.$$

In what follows we provide precise control of the constants $v_0$ and $b_0$ involved in inequality (A.4). It is easy to show that

$$v_0 \leq \frac{\|Kf\|_\infty}{t} \sup_{\{a \in \mathbb{R}^{jd}; \|a\|_2 < 1\}} \int \left(\sum_{|\lambda|<j} a_\lambda u_\lambda^j\right)^2 d\mu.$$

For any vector $a \in \mathbb{R}^{jd}$ we have

$$(A.5) \qquad \int \left(\sum_{|\lambda|<j} a_\lambda u_\lambda^j\right)^2 d\mu = \sum_{|\lambda|<j, |\lambda'|<j} a_\lambda a_{\lambda'} \int u_\lambda^j u_{\lambda'}^j \, d\mu$$

$$(A.6) \qquad\qquad\qquad \leq \sum_{|\lambda|<j} a_\lambda a_{\lambda'} \|u_\lambda^j\|_{L^2} \|u_{\lambda'}^j\|_{L^2}.$$

As argued in [5], the ellipticity property (2.1) yields

$$\|u_\lambda^j\|_{H^{-\nu/2}}^2 \leq \langle Ku_\lambda^j, u_\lambda^j \rangle = \langle \psi_\lambda, u_\lambda^j \rangle$$

$$\leq \|\psi_\lambda\|_{L^2} \|u_\lambda^j\|_{L^2} = \|u_\lambda^j\|_{L^2},$$

and the inverse inequality (see [5]) states that $\|u_\lambda^j\|_{L^2} \leq 2^{\nu j/2} \|u_\lambda^j\|_{H^{-\nu/2}}$, which implies that (dividing by $\|u_\lambda^j\|_{L^2}$)

$$(A.7) \qquad\qquad\qquad \|u_\lambda^j\|_{L^2} \leq 2^{\nu j}.$$

Using the above bound in the inequality (A.6) we finally obtain

$$\int \left(\sum_{|\lambda|<j} a_\lambda u_\lambda^j\right)^2 d\mu \leq 2^{2\nu j} \left(\sum_{|\lambda|<j} a_\lambda\right)^2 \leq 2^{j(2\nu+d)} \|a\|_2^2.$$

It follows that

$$v_0 \le \|Kf\|_\infty \frac{2^{j(2\nu+d)}}{t}.$$

Now, by definition of the constants $A_j$ and since $\sum_{|\lambda|<j} a_\lambda u_\lambda^j \in V_j$, we have

$$\left\| \sum_{|\lambda|<j} a_\lambda u_\lambda^j \right\|_\infty \le A_j \left\| \sum_{|\lambda|<j} a_\lambda u_\lambda^j \right\|_{L^2},$$

and it follows that, for $\|a\|_2 \le 1$,

$$\left\| \sum_{|\lambda|<j} a_\lambda u_\lambda^j \right\|_{L^2} \le 2^{j(\nu+d)},$$

which combined with the statement of Lemma A.5 gives

$$b_0 = O\left( \frac{2^{j(\nu+(3/2)d)}}{t} \right).$$

Now, recall that $\mathrm{Var}(\frac{1}{t} \int u_\lambda^j \, dG) = \frac{1}{t} \int (u_\lambda^j)^2 Kf \, d\mu$. Hence, using again the bound in inequality (A.7) we have

$$(A.8) \qquad \sum_{|\lambda|<j} \mathrm{Var}\left( \frac{1}{t} \int u_\lambda^j \, dG \right) \le \frac{1}{t} \|Kf\|_\infty 2^{j(d+2\nu)},$$

which implies that

$$\mathbb{E}(\chi^2(S_j)) \le \frac{1}{t} \|Kf\|_\infty 2^{j(d+2\nu)}.$$

Combining (A.8) and (A.3) yields finally

$$\mathbb{E}\left( \sum_{|\lambda|<j} \left( \frac{1}{t} \int u_\lambda^j \, dG - \langle f, \psi_\lambda \rangle \right)^2 \right) \le 2\left( \frac{1}{t} \|Kf\|_\infty 2^{j(d+2\nu)} + \|f\|_{L^2}^2 2^{-2js} \right).$$

From the Cauchy–Schwarz inequality and expression (A.4) with $\varepsilon = w$ it follows that there exists a constant $C > 0$, such that, for any $w > 0$

$$\mathbb{P}\{\chi^2(S_j) \ge C(1+w)^2 (2^{j(\nu+d/2)}/\sqrt{t} + 2^{j(\nu+(3/2)d)}/t)^2\} \le \exp(-w).$$

Combining the above inequalities, and using the fact that $f$ is bounded in $L^2$, we get finally

$$\mathbb{P}\{\|\hat{\alpha}_t - \alpha_0\|_2^2 \ge C(1+w)^2 \rho_{j,t}\} \le \exp(-w).$$

It remains to set $\eta = (1+w)$ and recall that $\rho_{j,t} = (2^{j(\nu+d/2)}/\sqrt{t} + 2^{j(\nu+(3/2)d)}/t)^2 + 2^{-2js}$ to get

$$\mathbb{P}\{\|\hat{\alpha}_t - \alpha_0\|_2^2 \ge C_1 \eta^2 \rho_{j,t}\} \le \exp(-\eta).$$

Now, applying Lemma 3.2 with $\theta_0 = \theta_j^*$, $\alpha = \hat{\alpha}_t$ and $b = e^{\|\log(f_{j,\theta_j^*})\|_\infty}$, we have that $\|\log(f_{j,\theta_j^*}/\exp(P_j g))\|_\infty \leq \varepsilon_j$, and so $b \leq M_1 e^{\varepsilon_j + \gamma_j}$. Hence, if $\eta^2 \rho_j^t \leq \frac{1}{4e^2 b^2 A_j^2}$, that is, if $\delta_j^t \leq 1/\eta^2$, then except in the set above, Lemma 3.2 implies that $\hat{\theta}_t = \theta(\hat{\alpha}_t) \in \mathbb{R}^{2^{jd}}$ exists and satisfies

$$\Delta(f_{j,\theta_j^*}; f_{j,\hat{\theta}_t}) \leq 2eb\|\alpha - \alpha_0\|_2^2$$

$$\leq 2M_1 e^{1+\varepsilon_j+\gamma_j} \eta^2 \rho_{j,t},$$

which completes the proof. □

PROOF OF THEOREM 4.2. From the bounds for $A_j$, $D_j$ and $\gamma_j$ given by Lemma A.5 and since $s > d/2$, we obtain that $\gamma_{j(t)} \to 0$ as $t \to \infty$ and so $\varepsilon_{j(t)} = O(A_j \Delta_j) = O(2^{-j(t)(s-d/2)})$. Hence, $\varepsilon_{j(t)} \to 0$ as $t \to \infty$ which implies that $\delta_{j(t)}^t = O(2^{-j(t)(2s-d)})$. Since $\varepsilon_{j(t)} \to 0$ and $\delta_{j(t)}^t \to 0$ as $t \to \infty$, Theorem 4.1 implies that

$$\Delta(f; f_{j(t),\theta_{j(t)}^*}) \leq O(2^{-2j(t)s}),$$

while for the estimation error, we have that as $t \to \infty$, then with probability tending to 1, $f_{j(t),\hat{\theta}_t}$ exists and by the Borel–Cantelli lemma satisfies

$$\Delta(f_{j(t),\theta_{j(t)}^*}; f_{j(t),\hat{\theta}_t}) \leq O(2^{-2j(t)s}).$$

The result now follows from the Pythagorean-like relationship (Lemma 3.1)

$$\Delta(f; f_{j(t),\hat{\theta}_t}) = \Delta(f; f_{j(t),\theta_{j(t)}^*}) + \Delta(f_{j(t),\theta_{j(t)}^*}; f_{j(t),\hat{\theta}_t}). \qquad \Box$$

PROOF OF THEOREM 5.1. Note that

$$\|\delta_t(\hat{\alpha}_t) - \alpha_0\|_2^2 = \sum_{|\lambda|<j} (\delta_t(\hat{\alpha}_{t,\lambda}) - \langle f, \psi_\lambda \rangle)^2$$

$$\leq 2\left( \sum_{|\lambda|<j} (\delta_t(\hat{\alpha}_{t,\lambda}) - \langle u_\lambda^j, Kf \rangle)^2 + (\langle u_\lambda^j, Kf \rangle - \langle f, \psi_\lambda \rangle)^2 \right).$$

From the proof of Theorem 4.1 [see (A.3)], we have (with the assumed conditions on $\psi$ and $K$)

$$\sum_{|\lambda|<j} (\langle u_\lambda^j, Kf \rangle - \langle f, \psi_\lambda \rangle)^2 \leq \|f\|_{L^2}^2 2^{-2js}.$$

Note that the space $B_{p,p}^s$ is continuously embedded in $H^\zeta$ whenever $\zeta \leq s + d/2 - d/p = 2\nu s/(2\nu + d)$. Moreover, since $2\nu/(2\nu + d) < 1$ and $f$ is uniformly bounded, we therefore obtain the estimate

$$\sum_{|\lambda|<j} (\langle u_\lambda^j, Kf \rangle - \langle f, \psi_\lambda \rangle)^2 \leq O(2^{-4js\nu/(2\nu+d)}).$$

This gives the optimal order $(\frac{1}{t})^{2s/(2s+2\nu+d)}$ provided $j(t)$ is large enough so that $2^{-j(t)} \leq (\frac{1}{t})^{(\nu+d/2)/(\nu(2s+2\nu+d))}$. For $t > 1$, we have $(\frac{1}{t})^{1/(2\nu)} \leq (\frac{1}{t})^{(\nu+d/2)/(\nu(2s+2\nu+d))}$, since $s \geq 0$, with equality if $s = 0$. Therefore, if $2^{-j(t)} \leq (\frac{1}{t})^{1/(2\nu)}$, we obtain

$$(A.9) \qquad \sum_{|\lambda| < j} (\langle u_\lambda^j, Kf \rangle - \langle f, \psi_\lambda \rangle)^2 \leq O\left(\left(\frac{1}{t}\right)^{2s/(2s+2\nu+d)}\right).$$

Note also that

$$\sum_{|\lambda| < j} (\delta_t(\hat{\alpha}_{t,\lambda}) - \langle u_\lambda^j, Kf \rangle)^2$$

$$= \sum_{|\lambda| < j} \left[ \Psi_\lambda^T K_j^{-1} \left( T_{\varepsilon(t)} \left( \frac{1}{t} \int \psi_\mu \, dG \right) - \langle h, \psi_\mu \rangle \right)_{|\mu| < j} \right]^2$$

$$= \left\| K_j^{-1} \left( T_{\varepsilon(t)} \left( \frac{1}{t} \int \psi_\mu \, dG \right) - \langle h, \psi_\mu \rangle \right)_{|\mu| < j} \right\|_2^2,$$

where $h = Kf$. Since $K$ is an isomorphism between $L^2$ and $H^\nu$, using the proof of Theorem 1 of [5] we have that for any $U = (u_\lambda)_{|\lambda| < j} \in \mathbb{R}^{2jd}$,

$$(A.10) \qquad \|K_j^{-1} U\|_2^2 \leq C \|U\|_{H^\nu} = \sum_{|\lambda| < j} 2^{2\nu|\lambda|} |u_\lambda|^2.$$

Hence, it follows that

$$\sum_{|\lambda| < j} (\delta_t(\hat{\alpha}_{t,\lambda}) - \langle u_\lambda^j, Kf \rangle)^2 \leq \sum_{|\lambda| < j} 2^{2\nu|\lambda|} \left( T_{\varepsilon(t)} \left( \frac{1}{t} \int \psi_\lambda \, dG \right) - \langle h, \psi_\lambda \rangle \right)^2.$$

We remark that $T_{\varepsilon(t)}(\frac{1}{t} \int \psi_\lambda \, dG) - \langle h, \psi_\lambda \rangle$ is exactly the error when estimating $\langle h, \psi_\lambda \rangle$ by the thresholding procedure on the "data" $\frac{1}{t} \int \psi_\lambda \, dG$. We are thus left with finding a threshold $\varepsilon(t)$ and some appropriate bounds for the estimation of $h$ based on the thresholded coefficients $T_{\varepsilon(t)}(\frac{1}{t} \int \psi_\lambda \, dG)$, $|\lambda| < j$, which would yield a bound for $\mathbb{E}\|\delta_t(\hat{\alpha}_t) - \alpha_0\|_2^2$.

To simplify the notation set $\hat{\beta}_{\lambda,t} = \frac{1}{t} \int \psi_\lambda \, dG$ and let $\beta_{\lambda,0} = \int \psi_\lambda h \, d\mu$. Since $G$ is a Poisson process with intensity $th$ it is easy to see that $\mathbb{E}(\hat{\beta}_{\lambda,t}) = \beta_{\lambda,0}$ and that $\text{Var}(\hat{\beta}_{\lambda,t}) = \frac{1}{t} \int \psi_\lambda^2 h \, d\mu = \frac{1}{t} \sigma_\lambda^2$. In order to bound the intensity estimation risk by a corresponding white-noise model risk, we will apply Lemma V of [3] to construct an approximation $\hat{\eta}_{\lambda,t}$ having an *exact* Gaussian distribution with the same mean $\beta_{\lambda,0}$ and the same variance $\text{Var}(\hat{\beta}_{\lambda,t})$. To this end, let $g_\lambda = \frac{1}{\sigma_\lambda} \psi_\lambda$ and note that $\int g_\lambda^2 h \, d\mu = 1$ and $\|g_\lambda\|_\infty = \frac{1}{\sigma_\lambda} \|\psi_\lambda\|_\infty = \frac{2^{|\lambda|d/2}}{\sigma_\lambda} := H_\lambda$, say. We construct $\hat{\eta}_{\lambda,t} = \beta_{\lambda,0} + t^{-1/2} \sigma_\lambda Z_\lambda$ by the following recipe.

First, if $\sigma_\lambda \geq C2^{|\lambda|d}\frac{\log^3 t}{t}$, then use Lemma V of [3] to construct $Z_\lambda$ and note that

$$V_{\lambda,t} = \mathbb{E}(\hat{\beta}_{\lambda,t} - \hat{\eta}_{\lambda,t})^2 = \frac{\sigma_\lambda^2}{t}\mathbb{E}(t^{-1/2}S_{\lambda,t} - Z_\lambda)^2 \leq C2^{|\lambda|d}t^{-2},$$

where $S_{\lambda,t} = \int \psi_\lambda(dG - th\,d\mu)$.

Second, if $\sigma_\lambda < C2^{|\lambda|d}\frac{\log^3 t}{t}$, choose an independent $Z_\lambda \sim N(0,1)$ and simply use the inequality

$$V_\lambda \leq 2\operatorname{Var}(\hat{\beta}_{\lambda,t}) + 2t^{-1}\sigma_\lambda^2 \leq C2^{|\lambda|d}\frac{\log^3 t}{t^2}.$$

In either case, we have therefore for all $|\lambda| < j$ and all $t > 0$,

$$V_{\lambda,t} \leq C2^{|\lambda|d}\frac{\log^3 t}{t^2}.$$

To apply the Gaussian approximation to $T_{\varepsilon(t)}(\hat{\beta}_{\lambda,t})$ note that

$$\mathbb{E}(T_{\varepsilon(t)}(\hat{\beta}_{\lambda,t}) - \beta_{\lambda,0})^2 \leq 2\mathbb{E}(T_{\varepsilon(t)}(\hat{\beta}_{\lambda,t}) - T_{\varepsilon(t)}(\hat{\eta}_{\lambda,t}))^2 + 2\mathbb{E}(T_{\varepsilon(t)}(\hat{\eta}_{\lambda,t}) - \beta_{\lambda,0})^2.$$

Since the mapping $y \to T(y, \varepsilon)$ is a contraction (see [8]) regardless of the value of $\varepsilon$, it follows that

$$\mathbb{E}(T_{\varepsilon(t)}(\hat{\beta}_{\lambda,t}) - \beta_{\lambda,0})^2 \leq 2V_{\lambda,t} + 2r(\varepsilon(t); t^{-1/2}\sigma_\lambda; \beta_{\lambda,0}),$$

where $r(\varepsilon(t); \sigma; \beta)$ is the Gaussian mean squared error $\mathbb{E}(T_\varepsilon(\beta + \sigma Z) - \beta)^2$ for estimation of $\beta$ from a single Gaussian observation with mean $\beta$ and variance $\sigma^2$. Since all intensities $h \in F_{p,p}^s(M)$ are uniformly bounded, we have $\sigma_\lambda \leq \|h\|_\infty$ and therefore

$$(\text{A.11}) \qquad \mathbb{E}(T_{\varepsilon(t)}(\hat{\beta}_{\lambda,t}) - \beta_{\lambda,0})^2 \leq 2V_{\lambda,t} + 2r(\varepsilon(t); t^{-1/2}\|h\|_\infty; \beta_{\lambda,0}).$$

Using the level-dependent threshold $\varepsilon(t) = 2^{\nu|\lambda|}t^{-1/2}\sqrt{|\log t|}$, the upper bound in inequality (A.11), the fact that $h$ belongs to a Besov ball $B_{p,p}^{s+\nu}(\tilde{M})$ for some finite constant $\tilde{M}$ and the stability property (A.10), we obtain the rate

$$\mathbb{E}\left(\sum_{|\lambda| < j}(\delta_t(\hat{\alpha}_{t,\lambda}) - \langle u_\lambda^j, Kf\rangle)^2\right) = O\left(\left(\frac{1}{t}\sqrt{|\log t|}\right)^{2s/(2s+2\nu+d)}\right),$$

as a particular case of classical results on soft wavelet thresholding (e.g., see [4]). Combining the above upper bound with inequality (A.9) and using again Lemma 3.2 concludes the proof of the theorem. □

**Acknowledgments.** We thank Marc Hoffmann for interesting and stimulating discussions and also for providing us the Matlab numerical procedures used in the paper by Cohen, Hoffmann and Reiss [5]. Both authors are very much indebted to the two anonymous referees, the Associate Editor and the Co-Editor Jianqing Fan, for their constructive criticism, comments and remarks that resulted in a major revision of the original manuscript.

LABORATOIRE IMAG-LMC                         DEPARTMENT OF STATISTICS
UNIVERSITY JOSEPH FOURIER                     UNIVERSITY PAUL SABATIER
BP 53                                         TOULOUSE
38041 GRENOBLE CEDEX 9                        FRANCE
FRANCE                                        E-MAIL: Jeremie.Bigot@cict.fr
E-MAIL: Anestis.Antoniadis@imag.fr